\documentclass[a4paper,12pt,draft]{article}
\usepackage{amsmath,amsfonts,amssymb,amsthm}

\newcommand\R{\mathbb{R}}
\newcommand\Z{\mathbb{Z}}

\def\dist{\mathrm{dist}\,}
\def\ep{\varepsilon}
\def\de{\delta}
\def\Lip{\mathrm{Lip}}
\def\diam{\mathrm{diam}\,}
\def\supp{\mathrm{supp}\,}

\begin{document}

\begin{center}
{\large 
\textbf{INVERSE SHADOWING AND RELATED MEASURES}}

\bigskip

{Sergey Kryzhevich\footnote{Corresponding author, E-mail:kryzhevicz@gmail.com}, Sergey Pilyugin
\bigskip 

St. Petersburg State University, Russia}
\end{center}

\noindent{\footnotesize
\textbf{Abstract.}
We study various weaker forms of inverse shadowing property for discrete dynamical systems on a smooth compact manifold. First, we introduce the so-called Ergodic Inverse Shadowing property (Birhhoff averages of continuous functions along the exact trajectory and the approximating one are close). We demonstrate that this property implies continuity of the set of invariant measures in Hausdorff metrics. We show that the class of systems with Ergodic Inverse Shadowing is quite broad, it includes all diffeomorphisms with hyperbolic nonwandering sets. Secondly, we study the so-called Individual Inverse Shadowing (any exact trajectory can be traced by approximate ones but this shadowing is not uniform with respect to selection of the initial point of the trajectory). We demonstrate that this property is closely related to Structural Stability and $\Omega$-stability of diffeomorphisms.

\bigskip

}

\bigskip

{\bf 1. Introduction}
\bigskip

The theory of shadowing of approximate trajectories (pseudotrajectories) of dynamical systems is now a well-developed field of the global theory of dynamical systems. Let us refer to the monographs [1--3] concerning the basics of the modern shadowing theory.

In parallel to the study of various (direct) shadowing properties, the so-called inverse shadowing properties were introduced (see [4,5]) and studied (see, for example, [6,7]). 

Recently, several authors studied the sets of shadowable points of dynamical systems in the context of their metric properties (see [8,9]).

Classical shadowing/inverse shadowing properties are closely related to Structural Stability and there are many interesting examples of systems without shadowing or inverse shadowing. Here we introduce weaker forms of inverse shadowing:
\begin{itemize}
\item we study inverse shadowing "almost always"\ -- the so-called Ergodic Inverse Shadowing property; this idea was inspired by the approach of the paper [10] where the so-called Ergodic Shadowing was introduced;
\item we introduce a "non-uniform" version inverse shadowing, the so-called Individual Inverse Shadowing. 
\end{itemize}

We study several properties of measures related to introduced forms of shadowing.

The structure of the paper is as follows. We give necessary definitions and formulate our main results in Sec.~2. Section~3 is devoted to Ergodic Inverse Shadowing property (EIS). We describe the class of systems with EIS.  Particularly, we demonstrate that this class is quite broad, it contains all systems with hyperbolic nonwandering sets. On the other hand, this property implies continuous dependence of the set of invariant measures with respect to small perturbations of the system (for this purpose, we spread the concept of invariant measures to non-autonomous discrete systems i.e. methods). In Sec.~4, we study metric properties of dynamical systems related to the so-called individual inverse shadowing property. We demonstrate that the $C^1$ interior of the introduced class of maps coincides with structurally stable diffeomorphisms. Also, we study a broader class of maps where all invariant measures are \emph{compatible} with inverse shadowing (Definition 4.1). This a weaker form of shadowing "almost everywhere". The $C^1$ -- interior of the latter set of systems coincides with the set of $\Omega$ -- stable diffeomorphisms. 
\bigskip

{\bf 2. Definitions and main results}
\medskip

Let $(K,\dist)$ be a compact metric space.

In this paper, we work with semi-dynamical systems (SDS) generated by surjective continuous mappings $f:\,K\to K$ and with dynamical systems (DS) generated by homeomorphisms $f:\,K\to K$.

Let $C(K,K)$ be the space of continuous mappings $g:\,K\to K$ with the metric
$$\dist_{C^0}(g,h)=\max_{x\in K}\dist(g(x),h(x)).$$

By definition, a {\em method} is a sequence
$$g=\{g_k\in C(K,K)\}, $$
where $k\in\Z_+$ in the case of SDS and $k\in\Z$ in the case of DS.

A sequence $\{y_k\}$ is called a {\em trajectory} of a method $g=\{g_k\}$ if
$$y_{k+1}=g_k(y_k),$$
where $k\in\Z_+$ in the case of SDS and $k\in\Z$ in the case of DS.

Fix a $d>0$. We say that a sequence $g=\{g_k\}$ is a $d$-{\em method} for a mapping $f:\,K\to K$ if
$$\dist_{C^0}(g_k,f)\leq d,$$
where $k\in\Z_+$ in the case of SDS and $k\in\Z$ in the case of DS.
\medskip

There are several different definitions of the inverse shadowing property. We will work with the definition used in [6] (there the methods which we use in this paper were called methods of the class $\Theta_s$).

Let us give this definition in the case of a DS $f:\,K\to K$.

We say that a homeomorphism $f$ has the {\em inverse shadowing property} if for any $\ep>0$ there exists a $d=d(\ep)>0$ such that for any trajectory $\{x_k=f^k(x):\,k\in\Z\}$ of $f$ and for any $d$-method  $g=\{g_k\}$ for $f$ there exists a trajectory $\{y_k\}$ of $g$ such that
$$\dist(x_k,y_k)\leq\ep,\quad k\in\Z. \eqno (2.1)$$

We denote by IS the set of homeomorphisms $f:\,K\to K$ having the inverse shadowing property.
\medskip

Our main goal is to study several modifications of the inverse shadowing property and relate them to properties of measures on the space $K$.
\medskip

The first of these modifications is related to the so-called ergodic approach to shadowing. Let us mention the paper [10] in which the authors introduced and studied the notionof ergodic (direct) shadowing. 
\medskip

We define a new property -- ergodic inverse shadowing property. We define it for SDS generated by surjective continuous mappings $f:\,K\to K$.

Denote by $\Lip_1$ the set of Lipschitz continuous functions $\varphi:\,K\to\R$ whose Lipschitz constant does not exceed 1.
\medskip

We say that a mapping $f$ has the {\em  ergodic inverse shadowing property} (EIS) if for any $\ep>0$ there exists a $d=d(\ep)>0$ such that for any trajectory $\{x_k=f^k(x):\,k\in\Z_+\}$ of $f$ and for any $d$-method $g$ for $f$ there exists a trajectory $\{y_k:\,k\in\Z_+\}$ of $g$ such that  
$$ \limsup_{n\to \infty} \dfrac1{n}  \sum_{k=1}^n (\varphi(x_k)-\varphi(y_k))\leq \ep \eqno (2.2)$$
for any function $\varphi\in\Lip_1$.
\medskip

We denote by EIS the set of mappings $f\in C(K,K)$ having the ergodic inverse shadowing property.
\medskip

In Sec. 3, we introduce the set of invariant measures for a method $g=\{g_k\}$. We introduce the notion of weak continuity of the set of Borel probability invariant measures of a mapping $f\in C(K,K)$ and show that if $f\in\mbox{EIS}$, then its set of probability  invariant measures is weakly continuous (Theorem~3.1).

In addition, we study some classical properties  of topological dynamics for mappings $f\in\mbox{EIS}$. For example, we show that if $f\in\mbox{EIS}$, then for any $\ep>0$ there exists a $d>0$ such that if $g\in C(K,K)$ and $\dist_{C^0}(f,g)<d$, then the $\ep$-neighbourhood of any minimal point of $f$ contains a Poisson stable point of $g$ (Proposition~3.4).

We also note that if the nonwandering set of a diffeomorphism $f$ of a closed smooth manifold is hyperbolic, then  $f\in\mbox{EIS}$ (Proposition~3.5).
\medskip

One of the principal differences between the (direct) shadowing property and inverse shadowing property is as follows: It is senseless to study the shadowing property selecting a single pseudotrajectory while, selecting an exact trajectory of a system, it is natural to study the inverse shadowing property for this selected trajectory.

Our second approach in this paper is based on the study of dynamical systems whose trajectories have the individual inverse shadowing property. 

The main definition is as follows. Let $f$ be a homeomorphism of $K$ and let $Y\subset K$.

We say that $f$ has the {\em individual inverse shadowing property} on the set $Y$ (and write $f\in \mbox{IIS(Y)}$) if  for any $x\in Y$ and $\ep>0$ there exists a $d>0$ such that for any $d$-method $g$ for which there exists a a trajectory $\{y_k:\,k\in\Z\}$ of $g$ such that inequalities (2.1) hold.

Note that in this case, $d$ depends not only on $\ep$ but also on the point $x$. 

If $Y=K$, we say that  $f$ has the individual inverse shadowing property (IIS) (and write $f\in \mbox{IIS}$).

In Sec. 4, we show that that the $C^1$-interior  of the set of diffeomorphisms  of a closed smooth manifold having the IIS coincides with the set of structurally stable diffeomorphisms (Theorem~4.1). This result generalises the main statement of the paper [11] concerning diffeomorphisms having the inverse shadowing property.

In the same Sec.~4, we study relations between individual inverse shadowing and measures on the space $K$.

Let ${\mathcal M}(K)$ be the set of all nonatomic probability Borel measures on the space $K$.

Introduce the following notation: for $\ep,d>0$ and a  homeomorphism $f$ of $K$, let $\Phi(\ep,d,f)$ be the set of all points $x\in K$ such that for any $d$-method $g$ for $f$ there exists a trajectory $\{y_k\}$ of $g$ that satisfies inequalities (2.1).

We say that a measure $\mu\in{\mathcal M}(X)$ is {\em compatible with inverse shadowing for $f$} if for any $\ep>0$  there exists a $d>0$ such that if $\mu(A)>0$, then
$$ A\cap\Phi(\ep,d,f)\neq\emptyset.$$

We prove the following two statements:

$\bullet$ If a strictly positive measure $\mu\in{\mathcal M}(K)$ is compatible with inverse shadowing for $f$, then $f\in \mbox{IS}$ (Corollary~4.1);

$\bullet$ the $C^1$-interior of the set of diffeomorphisms $f$ of a closed smooth manifold for which there exists a strictly positive measure compatible with inverse shadowing for $f$ coincides with the set of structurally stable diffeomorphisms (Corollary~4.2).

Finally, for a homeomorphism $f$ of the space $K$, we define the set
$$\Psi(f)=\{x\in K:\, f\in\mbox{IIS}(\{x\})\}.$$ 

Thus, $f$ has the individual inverse shadowing property on a set $Y\subset K$ if and only if $Y\subset\Psi(f)$.

We say that a measure $\mu\in{\mathcal M}(K)$ is {\em compatible with individual inverse shadowing for $f$} if for any set $A\subset K$  with $\mu(A)>0$, 
$$A\cap\Psi(f)\neq\emptyset.$$

We show that the $C^1$-interior of the set of diffeomorphisms $f$ of a closed smooth manifold $M$ for which every $f$-invariant measure $\mu\in{\mathcal M}(M)$ is compatible with individual inverse shadowing coincides with the set of $\Omega$-stable diffeomorphisms (Theorem~4.2).
\bigskip

{\bf 3. Ergodic inverse shadowing and invariant measures}
\bigskip

Let $\mathcal M$ be the set of all Borel probability measures on $K$ with the $*$-weak convergence topology. Introduce the so-called Kantorovich--Wasserstein metric $\rho$ on $\mathcal M$ as follows:
$$ W_1(\mu_1,\mu_2)=\sup_{\varphi\in \mathrm{Lip}_1} \left| \int_K \varphi\, d\mu_1 - \int_K \varphi\, d\mu_2\right| $$ 
(recall that $\Lip_1$ is the set of Lipschitz continuous functions $\varphi:\,K\to\R$ whose Lipschitz constant does not exceed 1).The convergence in the Kantorovich--Wasserstein metric is equivalent to the  $*$-weak convergence [12]. 

Let $C(K,\R)$ be the space of continuous functions on $K$ with the metric
$$\dist_{C^0}(f,g)=\max_{x\in K}|f(x)-g(x)|.$$
 
Given a continuous mapping $f:K\to K$, we define the so-called push-forward map $f^\#: \mathcal M \to  \mathcal M$ as follows: $f^\#(\mu)=\nu$   if 
$$\int_{K} \varphi \circ f\, d\mu=\int_{K} \varphi \, d\nu$$ 
for any $\varphi \in C(K,{\mathbb R})$. The operator $f^\#$ is linear continuous, and $f^\# \delta(x)=\delta (f(x))$ for any $x\in K$, where $\delta(x)$ is the Dirac measure taken at the point $x$.  Recall that a measure $\mu$ is invariant with respect to $f$  if and only if $f^\#(\mu)=\mu$.

For $\varepsilon>0$ and $A\subset{\mathcal M}$, we denote by $U_\varepsilon (A)$ the $\varepsilon$-neighbourhood  of $A$ in the metric $\rho$.

First of all, we introduce the set of invariant measures for a method.

Let $g=\{g_k:\,k\in\Z_+\}$ be a method; denote  $g^k_0=g_{k-1}\circ \ldots \circ g_0$.

For a point $x\in K$, we define the set
$${\mathcal M}_0 (g,x)= \bigcap_{n=1}^\infty \overline{ \left\{\dfrac1N \sum_{k=0}^{N-1} \delta(g^k_0(x)) \, : \,N\ge n\right\}}. $$
In other words, this is the set of limit points for the sequence 
$$\left\{\dfrac1N \sum_{k=0}^{N-1}  \delta(g^k_0(x)) : N\in {\mathbb N}\right\}. $$
Note that  that all the sets ${\mathcal M}_0 (g,x)$  are nonempty as intersections of nested nonempty compact sets. 

Let 
$${\mathcal M}_0(g):=\bigcup_{x\in K} {\mathcal M}_0 (g,x) \mbox{ and } {\mathcal M}(g):= \overline{\mathrm{conv}\, {\mathcal M}_0 (g)}.$$ 
We call any measure of the set ${\mathcal M}(g)$ \emph{invariant} with respect to the method $g$. We start with a folklore statement which is a corollary of the Birkhoff ergodic theorem. 
\medskip

{\bf Proposition~3.1. } {\em Let $f$ be a continuous mapping of a compact metric space $K$ into itself and let $\mu$ be an ergodic $f$-invariant probability measure. Then for $\mu$-almost all points $x$ of $K$, the sequence 
$$\dfrac1n \sum_{k=0}^{n-1} \delta (f^k(x))$$
converges to $\mu$ $*$-weakly}.
\medskip

{\em Proof.} Consider a countable set $\{\varphi_j: j\in {\mathbb N}\}$ of continuous functions that is dense in $C(K, {\mathbb R})$. By Birkhoff's ergodic theorem, there exists a set $X\subset K$ of full measure such that for any $x\in X$ and any $j\in {\mathbb N}$,
$$\dfrac1n \sum_{k=0}^{n-1} \varphi_j (f^k(x)) \to \int_K \varphi_j \, d\mu. $$
We demonstrate that a similar statement is true for any function 
$$\varphi\in C(K, {\mathbb R}).$$ 
Fix a point $x\in X$, a continuous function $\varphi$, an arbitrary $\varepsilon>0$, and take an index $j\in {\mathbb N}$ such that  $\dist_{C^0}(\varphi_j,\varphi)<\varepsilon$. Then 
$$\left| \int_M \varphi\, d\mu-  \int_M \varphi_j\, d\mu\right|<\varepsilon.$$
Also, there exists a number $N\in {\mathbb N}$ such that 
$$\left|\dfrac1n \sum_{k=0}^{n-1} \varphi (f^k(x))-  \dfrac1n \sum_{k=0}^{n-1} \varphi_j (f^k(x)) \right|\le \varepsilon$$
for any $n\ge N$.

These relations imply that
$$\begin{array}{c}
\int\limits_K \varphi \, d\mu -2\varepsilon\le \liminf  \dfrac1n \sum\limits_{k=0}^{n-1} \varphi (f^k(x)) \le\\
\limsup  \dfrac1n \sum\limits_{k=0}^{n-1} \varphi (f^k(x)) \le \int\limits_K \varphi \, d\mu +2\varepsilon,
\end{array}$$
from which our statement follows. $\Box$
\medskip

{\bf Proposition~3.2. }{\em If all the mappings  $g_k$ of the method $g$ are the same, i.e. they coincide with a given map $g_0$, then the set ${\mathcal M}(g)$ is the set of $g_0$-invariant measures in the classical sense.}
\medskip

{\em Proof.} First of all, we show that any measure in the set ${\mathcal M}(g)$ is invariant. Indeed, following the lines of the Krylov--Bogolyubov theorem ([13, Theorem~4.1.1]) we check that any measure in any set ${\mathcal M}_0 (g,x)$, $x\in K$, is $g_0$-invariant. In addition, the set of invariant measures is always closed in the $*$-weak topology and convex. Therefore, the set ${\mathcal M}(g)$ is a subset of the set of all $g_0$-invariant measures.

Conversely, if $\mu$ is an ergodic invariant measure of $g_0$, then, by the Birkhoff ergodic theorem, for $\mu$-almost all points $x$, the set ${\mathcal M}_0(g,x)$ is the singleton $\{\mu\}$. Thus, any ergodic $g_0$-invariant measure is an element of ${\mathcal M}(g)$. On the other hand, any invariant measure is an element of the convex hull of the set of all ergodic measures, and, therefore, it is an element of ${\mathcal M}(g)$. $\Box$
\medskip

We show that the set ${\mathcal M}(g)$ is upper semicontinuous in the Hausdorff metric in the following sense. 

For  $\varepsilon>0$ and  $A\subset{\mathcal M}$, we denote by $U_\varepsilon (A)$ the $\varepsilon$-neighbourhood of $A$ in the Kantorovich-Wasserstein metric $\rho$.
\medskip

{\bf Lemma~3.1. }{\em Consider a sequence of methods $g^{(k)}=\{g_{m}^{(k)}: k,m \in {\mathbb Z}_+\}$ for which  there exists a mapping $f$ such that $g_{m}^{(k)}\rightrightarrows  f$ as $k\to \infty$ uniformly with respect to $m\in {\mathbb Z}_+$. Then for any $\alpha>0$ there exists an $N\in {\mathbb Z}_+$ such that 
$${\mathcal M}(g^{(k)})\subset U_\alpha ({\mathcal M}(f))$$ 
for any} $k\ge N$.
\medskip

{\em Proof.} Take a sequence of methods $g^{(k)}$ satisfying the conditions of the lemma and consider a sequence $\mu_k$ of $g^{(k)}$-invariant measures  $*$-weakly converging to a measure $\mu^*$. Then the sequence of push-forward measures $f^\# \mu_k$ converges to the measure $f^\# \mu^*$. 

Let us demonstrate that the measures $\mu^*$ and $f^\# \mu^*$ coincide, and, therefore, $\mu^* \in {\mathcal M}(f)$. Fix an $\varepsilon>0$. Let a number $k$ be so large that 
$$W_1(\mu_k, \mu^*)<\varepsilon/8,\qquad W_1(f^\#\mu_k, f^\#\mu^*)<\varepsilon/8,$$ and $\dist_{C^0}(f,{g^{(k)}})<\varepsilon/8$ for all $m \in {\mathbb Z}_+$. 

Since $\mu_k$ is an element of the closure of the convex hull of ${\mathcal M}_0(g^{(k)})$, we can select a finite convex combination 
$${\widetilde{\mu_k}}:=\sum_{j=1}^n \alpha_j \mu_{k,j},$$
where $\mu_{k,j}\in {{\mathcal M}_0(g^{(k)})}$, all $\alpha_j>0$,  
$$\sum\limits_{j=1}^n \alpha_j =1,\quad W_1(\mu_k, \widetilde{\mu_k})<\varepsilon/8,\quad\mbox{and}\quad W_1(f^\# \mu_k, f^\# \widetilde{\mu_k})<\varepsilon/8.$$

In general, the number $n$ depends on the measure $\mu_k$ and on $\varepsilon$.

Thus, it suffices to prove that 
$$W_1(\mu_{k,j}, f^\# \mu_{k,j})<\varepsilon/2  \quad \mbox{for any}\quad j=1,\ldots,n.\eqno (3.1)$$

Indeed, if inequality (3.1) is true, then  $W_1(\mu^*,f^\# \mu^*)\le \varepsilon$, which  completes  the proof since $\varepsilon$ is arbitrarily small.

Let us check inequality (3.1). Take an index $j\in \{1,\ldots,n\}$ and the corresponding measure $\mu_{k,j}$. By definition, there exists a point $x \in K$ and a sequence $m_l \to \infty$ such that 
$$\mu_{k,j}=  \lim_{l\to \infty} \mu_{k,j,l}, \quad \mbox{where} \quad \mu_{k,j,l} = \dfrac1{m_l} \sum_{i=0}^{m_l-1} \delta(g^i_0 (x)).$$

Take a number $L\in {\mathbb N}$ so large that 
$$W_1(\mu_{k,j}, \mu_{k,j,L})<\varepsilon/8, \quad W_1(f^\#\mu_{k,j}, f^\# \mu_{k,j,L})<\varepsilon/8,$$
and, also, $\diam K/m_L <\varepsilon/8$. Then
$$W_1(\mu_{k,j,L},f^\#\mu_{k,j,L})= \dfrac1{m_L}W_1  \left(\sum\limits_{i=0}^{m_L-1} \delta ([g^{k}]_0^i (x)), \sum\limits_{i=0}^{m_L-1} \delta (f\circ [g^{k}]_0^i (x))\right)\le $$
$$ \le\dfrac1{m_L}W_1(\delta(x),\delta(g^{(k)}_{m_L} \circ [g^{(k)}]_0^{m_L-1} (x))) + $$
$$+\dfrac1{m_L}\sum\limits_{i=1}^{m_L-1} W_1(\delta ([g^{(k)}]_0^i (x)), \delta(f\circ [g^{(k)}]_0^{i-1}(x))) < \dfrac{\varepsilon}4.$$

This proves (3.1). 

Now, let the statement of the lemma be wrong. Then we may assume that there exists an $\alpha>0$ such that $\mu_k\notin U_\alpha ({\mathcal M}(f))$ for any $k\in {\mathbb N}$ (here we pass to a subsequence if necessary). We may assume, without loss of generality, that the sequence $\mu_k$ $*$-weakly converges to a measure $\mu^*$ (which is evidently outside $U_\alpha ({\mathcal M}(f))$). This contradicts to what have been proved above. $\Box$
\medskip

The ergodic inverse shadowing yields the converse inclusion -- the set of invariant measures of a mapping having the EIS property belongs to a small neighbourhood of the set of invariant measures for close methods. Namely, the following result is true. 
\medskip

Let us define the corresponding property.

We say that the set of Borel probability invariant measures of a continuous surjective mapping $f$ is {\em weakly continuous} if for any $\varepsilon>0$ there exists a$d>0$ such that if  $g=\{g_k\}$ is a $d$-method for $f$, then ${\mathcal M}(f)\subset U_\varepsilon ({\mathcal M}(g))$.

Let $\mathrm{CM}$ be the class of all mappings with weakly continuous sets of Borel probability invariant measures.
\medskip

{\bf Theorem~3.1. }{\em  If}  $f\in\mbox{EIS}$, {\em then} $f\in\mbox{CM}$. 
\medskip

{\em Proof.}  For any $\varepsilon>0$ there is a $d>0$ such that for any point $x\in K$, any $d$-method $g$ for $f$ and any function $\varphi \in \Lip_1$ there exists a point $y\in K$ such that
$$\limsup \dfrac1n \sum_{k=0}^{n-1} (\varphi(g^k_0(y))-\varphi(f^k(x)) < \varepsilon.\eqno (3.2)$$

Fix an $\varepsilon>0$ and a corresponding $d>0$ and consider a $d$-method $g$.

Take any ergodic $f$-invariant measure $\mu$. Then, for $\mu$-almost all points $x$, we have the equality
$$\mu=\lim \dfrac{1}{n} \sum_{j=0}^{n-1} \delta(f^j(x)).$$
Fix such a point $x$ and select a point $y$ that corresponds to $x$ in the sense of (3.2). Take a subsequence $\{n_k\}$ such that the sequence
$$\dfrac{1}{n_k} \sum_{j=0}^{n_k-1} \delta(g^j_0(y))$$
converges $*$-weakly to a measure $\nu$. Then, by (3.2), the Kantorovich--Wasser\-stein distance between $\mu$ and $\nu$ does not exceed $\varepsilon$. Since the set of all finite convex combinations of ergodic invariant measures is dense in ${\mathcal M}(f)$, this completes the proof. $\Box$  
\medskip

Let us recall two classical definitions of topological dynamics.
\medskip

We say that a point $x\in K$ is {\em Poisson stable} for a mapping $f$ if there is a sequence $n_k\to \infty$ such that $f^{n_k}(x)\to x$. 
\medskip

We say that a point $x\in K$ is {\em minimal} for a mapping $f$ if the set 
$$\overline{O^+_f(x)}:=\overline{\{f^k(x):k\in {\mathbb Z}_+\}}$$
is minimal, i.e., $\overline{O^+_f(y)}=\overline{O^+_f(x)}$ for any $y\in \overline{O^+_f(x)}$.
\medskip

{\bf Proposition~3.3.} {\em Let $x$ be a minimal point of a continuous mapping $f$ of a compact metric space. Then there exists an $f$-invariant probability measure $\mu$ such that} $x\in \supp\, \mu$.
\medskip

{\em Proof. } Take an invariant probability measure $\mu$ for the mapping $f|_{\overline{O_f^+(x)}}$. The support of $\mu$ is a closed invariant subset of $\overline{O_f^+(x)}$; thus, it must coincide with $\overline{O_f^+(x)}$. $\Box$
\medskip

Since any minimal point of any continuous mapping always belongs to the support of an invariant measure (actually defined by that point) and, on the other hand, any point of the support of any invariant measure is Poisson stable, we can formulate the following statement.
\medskip

{\bf Proposition~3.4.} {\em Assume that} $f\in\mbox{EIS}$. {\em Then for any $\varepsilon>0$ there exists a $d>0$ such that if} $\dist_{C^0}(f,g)<d$, {\em then the $\varepsilon$-neighbourhood of any minimal point of $f$ contains a Poisson stable point of} $g$.
\medskip

Finally, we indicate a class of diffeomorphisms of a closed smooth manifold having the EIS.
\medskip

{\bf Proposition~3.5. }{\em If the nonwandering set of a diffeomorphism $f$ is hyperbolic, then} $f\in\mbox{EIS}$.
\medskip

{\em Proof. } Let the nonwandering set $\Omega(f)$ of a diffeomorphism $f$ be hyperbolic. It is well known (see [14], for example) that there exists a neighbourhood $U$ of $\Omega(f)$ such that any segment of a trajectory of $f$ belonging to $U$ is hyperbolic with the same hyperbolicity constants $C\geq 1$ and $\lambda\in(0,1)$.

Fix an arbitrary point $x\in M$. Since the trajectory of $x$ tends to $\Omega(f)$ as time grows, there exists a number $N$ such that $f^k(x)\in U$ for $k\geq N$. The set
$$F:=\{f^k(x):\,k\geq N\}$$
belongs to $U$; its hyperbolicity implies that on $F$ there exists a $(C,\lambda)$-structure and for any $\varepsilon>0$ there exists a $d=d(\varepsilon,C,\lambda)$ with the following property: If $g=\{g_k\}$ is a $d$-method for $f$, then there exists a sequence $\{y_k:\,k\geq N\}$ such that
$$y_{k+1}=g_k(y_k)\quad\mbox{and}\quad \dist(y_k,x_k)\leq\varepsilon/2,\quad k\geq N$$
(see [7]).

To complete the proof of our proposition, take $y_0\in (g_0^{N-1})^{-1}(y_N)$ and note that inequality (2.2) is obviously valid for any function $\varphi:\,M\to\R$ whose Lipschitz constant does not exceed 1. $\Box$
\medskip

{\bf Remark.} It follows from the last result that if $M$  is a closed smooth manifold, then the set $\mbox{EIS}(M)$ contains the set of $\Omega$-stable diffeomorphisms (and hence, the $C^1$-interior of the set $\mbox{EIS}(M)$ has this property). Thus, the set $\mbox{EIS}(M)$ is strictly larger than the set $\mbox{IS}(M)$ (whose $C^1$-interior coincides with the set of structurally stable diffeomorphisms, see [11]).
\bigskip

{\bf 4. Individual inverse shadowing}
\bigskip

In this section, we study sets of diffeomorphisms of a closed smooth manifold $M$ having the individual inverse  shadowing property on various subsets of the phase space.

Let us first formulate two basic results from the theory of structural stability.

For a set $A$ of diffeomorphisms, denote by $\mbox{Int}^1(A)$ its interior with respect to the $C^1$-topology.
\medskip

Denote by ${\cal H}(M)$ the set of diffeomorphisms of $M$ for which any periodic point is hyperbolic; let ${\cal F}(M)=\mbox{Int}^1({\cal H}(M))$.

Denote by $\mbox{KS}(M)$ the subset of ${\cal H}(M)$ consisting of diffeomorphisms for which stable and unstable manifolds of periodic points are transverse (diffeomorphisms $f\in\mbox{KS}(M)$ are called Kupka--Smale diffeomorphisms).
\medskip

{\bf Proposition~4.1. }

(1) {\em The set } ${\cal F}(M)$ {\em coincides with the set} $\Omega\mbox{S}(M)$ {\em of $\Omega$-stable diffeomorphisms of} $M$.

(2) {\em The set } $\mbox{Int}^1(\mbox{KS}(M))$ {\em coincides with the set} $\mbox{StS}(M)$ {\em of structurally stable diffeomorphisms of} $M$.
\medskip

The first statement of Proposition~4.1 is proved in [15], the second one is proved in [16] (see also the book [3]).
\medskip

\medskip

In what follows, we do not mention the manifold $M$ and write ${\cal F}$, $\mbox{KS}$, $\Omega\mbox{S}$, and $\mbox{StS}$ instead of ${\cal F}(M)$, $\mbox{KS}(M)$, $\Omega\mbox{S}(M)$, and $\mbox{StS}(M)$, respectively.
\medskip

Our first goal is to prove the following statement.
\medskip

{\bf Theorem~4.1. }{\em The set} $\mbox{Int}^1(\mbox{IIS})$ {\em coincides with the set of structurally stable diffeomorphisms.}
\medskip

We divide the proof into several steps.
\medskip

{\bf Lemma~4.1. }{\em Assume that a diffeomorphism $f$ belongs to} $\mbox{Int}^1(\mbox{IIS})$. {\em Then any periodic point of $f$ is hyperbolic}.
\medskip

{\em Proof.} To get a contradiction, assume that a diffeomorphism $f\in\mbox{Int}^1(\mbox{IIS})$ has a nonhyperbolic periodic point $p$. Obviously, $f\in\mbox{Int}^1(\mbox{IIS})$ if and only if $f^n\in\mbox{Int}^1(\mbox{IIS})$ for some (every) natural $n$; hence, without loss of generality, we may assume that $p$ is a fixed point of $f$.

We apply a standard construction finding a diffeomorphism $h$ that is $C^1$-close enough to $f$ (so that $h\in\mbox{IIS}$) and linear in a neighbourhood of $p$; such constructions are described in detail in the book [3], and we apply them  several times in this paper.

To simplify presentation, we only consider the case where the derivative $Df(p)$ has an eigenvalue 1; the remaining possible cases are left to the reader (details can also be found in the book [3]).

Then we can find a diffeomorphism $h\in\mbox{IIS}$ and a neighbourhood $U$ of $p$ with local coordinates $(u,v)$ such that
 
 -- $p$ is the origin in $U$;
 
 -- the coordinate $u$ is one-dimensional and the coordinate $v$ is $(m-1)$-dimensional, where $m$ is the dimension of $M$;

$$R:=\{|u|,|v|\leq a\}\subset U, \eqno (4.1)$$
where $a>0$;

-- in $U$, $h$ has the form
$$h(u,v)=(u,Av) \eqno (4.2)$$
with some matrix $A$.

Take an arbitrary $\ep\in(0,a)$ and an arbitrary $d>0$. Clearly, there exists a $\de>0$ and a mapping $g\in C(M,M)$ such that
$$g(u,v)=(u+\de,Av) \mbox{ in } U \eqno (4.3)$$
and
$$\dist(g(x),h(x))\leq d,\quad x\in M. \eqno (4.4)$$

Let $g_k=f$ for $k<0$ and $g_k=g$ for $k\geq 0$; then $\xi=\{g_k:\,k\in\Z\}$ is a $d$-method for $h$.

Thus, if $(u,v)\in U$ and $k\geq 0$, then $g_k(u,v)=(u+\de,Av)$.

Clearly, for any point $x_0\in M$ there exists an 
$n\geq 0$ such that 
$$
g_n\circ\cdots\circ g_0(x_0)\notin R,
$$
i.e., for any trajectory $\{x_k\}$ of $\xi$ there exists an $n\geq 0$ such that $$\dist(x_{n+1},p)>\ep,$$ which means, due to the arbitrariness of $d$, that $h\notin\mbox{IIS}$.

The obtained contradiction completes the proof. $\Box$
\medskip

Thus, we have shown that $\mbox{Int}^1(\mbox{IIS})
\subset{\cal H}$.

Since the set $\mbox{Int}^1(\mbox{IIS})$ is obviously $C^1$-open, it follows from item (1)  of Proposition~4.1 that any diffeomorphism in $\mbox{Int}^1(\mbox{IIS})$ is $\Omega$-stable.
\medskip

{\bf Remark. } In fact, we have used the assumption $f\in\mbox{Int}^1(\mbox{IIS}(\mbox{Per}(f))$; thus, the $C^1$-interior of the set of diffeomorphisms having the individual inverse shadowing property on the set of their periodic points consists of $\Omega$-stable diffeomorphisms.
\medskip

{\bf Lemma~4.2. }{\em Let $p$ and $q$ be periodic points of  a diffeomorphism $f$ belonging to} $\mbox{Int}^1(\mbox{IIS}$). {\em Then the unstable manifold $W^u(p)$ of $p$ and the stable manifold $W^s(q)$ of $q$ are transverse}.
\medskip

{\em Proof. } To get a contradiction, assume that there is a point $r$ at which $W^u(p)$ and $W^s(q)$ are nontransverse. By the previous lemma, $f$ is $\Omega$-stable; it follows that
$p$ and $q$ belong to different basic sets of $f$.

As in the previous lemma, we assume for simplicity that $p$ and $q$ are fixed points.

It is shown in [17] that we can find a diffeomorphism $h\in\mbox{IIS}$  and select a point $r$ of nontransverse intersection of $W^u(p,h)$ and $W^s(q,h)$ (here $W^u(p,h)$ and $W^s(q,h)$ are the unstable manifold of $p$ and the stable manifold of $q$ for the diffeomorphism $h$) such that the following statements hold:

-- there exists a neighbourhood  $V$ of $q$ in which $q$ is the origin and such that $h$ is linear in $V$;

-- the points $h^k(r),\,k\geq 0,$  belong to $V$;

-- if $\mbox{dim} W^u(p,h)=u$, then there exists an $u$-dimensional linear subspace $L$ (with respect to local coordinates in $V$) and an $u$-dimensional disk $C$ with the following properties:

(c1) $C$ is open in the inner topology of the affine space $r+L$, and its closure belongs to the intersection of $r+L$ with $V$;

(c2) $r\in C$;

(c3) for any $\ep>0$ small enough, $C$ contains an open $u$-dimensional disk $C_{\ep}$ such that if
$$\dist(h^k(x),h^k(r))\leq \ep,\quad k\leq 0, $$
for a point $x\in C$, then $x\in C_{\ep}$.

Clearly, in this case, we can identify $L$ with $T_rW^u(p,h)$, the tangent  space of $W^u(p,h)$ at the point $r$.

Denote by $S$ and $U$ the linear subspaces (in local coordinates of $V$) such that $S\cap V=W^s(q,h)\cap V$ and $U\cap V=W^u(q,h)\cap V$.

The nontransversality of $W^u(p,h)$ and $W^s(q,h)$ at $r$ means that $L+S\neq T_rM$. If $\pi$ is the projection in $T_rM$ to $U$ parallel to $S$, then the above condition means that
$$\pi L\neq U. \eqno (4.5)$$

It follows from (4.5) that there exists a one-dimensional subspace $Z$ of $U$ such that $\pi L\cap Z=\{0\}$. Let $z$ be a unit vector in $Z$.

Let us show that $h\notin\mbox{IIS}(\{r\})$. Assume the converse and take $\ep>0$ so small that the closure of the $\ep$-neighbourhood of the set $\{h^k(r):\,k\geq 0\}$ belongs to $V$ (this is possible due to (c2)) and property (c3) is satisfied. Find the corresponding $d>0$.

Clearly, there exists a $g\in C(M,M)$ such that
$$g(x)=h(x)+\de z,\quad x\in C,$$
where $\de>0$, and 
$$\dist(g(x),h(x))<d,\quad x\in M.$$

Then the family $\xi=\{g_k\}$ with $g_k=h$ for $k\neq 0$ and $g_0=g$ is a $d$-method for $h$.

It follows from (c2) that if $\{x_k\}$ is a trajectory of $\xi$ for which inequalities (2.1) (with $f$ replaced by $h$) hold, then $x_0\in C_{\ep}$, but then $\pi(g(x_0))\notin \pi L$, which obviously implies that the points $x_n$ leave $V$ (and the closed $\ep$-neighbourhood of the positive $h$-trajectory of $r$) as $n$ grows.

The obtained contradiction completes the proof. $\Box$
\medskip

It follows from Lemmas~4.1 and 4.2 that
$$\mbox{Int}^1(\mbox{IIS})\subset\mbox{KS}.$$

Hence,
$$\mbox{Int}^1(\mbox{IIS})\subset\mbox{Int}^1(\mbox{KS}),$$
and item (2) of Proposition~4.1 implies that
$$\mbox{Int}^1(\mbox{IIS})\subset\mbox{StS}.$$

The converse inclusion has been proved in [6]. Thus, Theorem~4.1
is proved.
\medskip

Let ${\mathcal M}(K)$ be the set of all nonatomic probability Borel measures on a compact metric space $(K,\dist)$. Denote by $\mbox{supp}(\mu)$ the support of a measure $\mu\in{\mathcal M}(X)$.
\medskip

Recall the notation and definition introduced in Sec.~2.

For $\ep,d>0$ and a homeomorphism $f$ of $K$, let $\Phi(\ep,d,f)$ be the set of all points $p\in K$ such that for any $d$-method $\xi$ for $f$ there exists a trajectory $\{x_k\}$ of $\xi$ that satisfies inequalities (2.1).

Clearly, for a set $Y\subset K$, $f\in \mbox{IS}(Y)$ if and only if for any $\ep>0$  there exists a $d>0$ such that $Y\subset\Phi(\ep,d,f)$.
\medskip

{\bf Definition~4.1. } We say that a measure $\mu\in{\mathcal M}(K)$ is {\em compatible with inverse shadowing for $f$} if for any $\ep>0$ there exists a $d>0$ such that if $\mu(A)>0$, then
$$A\cap\Phi(\ep,d,f)\neq\emptyset. \eqno (4.6)$$

{\bf Lemma~4.3. }{\em If $\mu$ is compatible with inverse shadowing for $f$, then} $f\in \mbox{IS}(\mbox{supp}(\mu))$.
\medskip

{\em Proof.} Take a point $p\in\mbox{supp}(\mu)$ and a natural number $n$. Fix an arbitrary $\ep>0$.

Since $f$ and $f^{-1}$ are uniformly continuous, there exists a neighbourhood $V_n$ of $p$ such that
$$\dist(f^k(p),f^k(q))\leq\ep/2,\quad q\in V_n,\,-n\leq k\leq n. \eqno (4.7)$$

Since $p\in\mbox{supp}(\mu)$, $\mu(V_n)>0$.

Take $d$ corresponding to $\ep/2$ in the above definition of compatibility with inverse shadowing and apply (4.6) to find a point 
$$q_n\in V_n\cap \Phi(\ep/2,d,f).$$

Let $\xi=\{g_k\}$ be a $d$-method for $f$ and let $\{x^n_k:\,k\in \Z\}$ be a trajectory of $\xi$ such that
$$\dist(x_k^n,f^k(q_n))\leq\ep/2,\quad k\in \Z. \eqno (4.8)$$

It follows from (4.7) and (4.8) that
$$\dist(x_k^n,f^k(p))\leq\ep,\quad -n\leq k\leq n.$$

Apply the diagonal process to find a sequence $\{x_k\}$ such that $x^n_k\to x_k$ as $n\to\infty$ for any fixed $k\in\Z$.

Clearly, $g_k(x_k)=x_{k+1}$ and 
$$ \dist(x_k,f^k(p))\leq\ep,\quad k\in\Z, $$
which completes the proof. $\Box$
\medskip

{\bf Corollary~4.1. }{\em If $\mu$ is compatible with inverse shadowing for $f$ and} $\mbox{supp}(\mu)=K$ {\em (i.e., the measure $\mu$ is strictly positive), then} $f\in \mbox{IS}$.
\medskip

Corollary~4.1 and Theorem~4.1 imply the following statement.

{\bf Corollary~4.2. }{\em Let ${\mathcal I}_1$ be the set of diffeomorphisms $f$ of a closed smooth manifold for which there exists a strictly positive measure compatible with inverse shadowing for $f$.

Then the set} $\mbox{Int}^1({\mathcal I}_1)$ {\em coincides with the set of structurally stable difeomorphisms.}
\medskip

{\bf Corollary~4.3. }{\em Let ${\mathcal I}_2$ be the set of diffeomorphisms $f$ of a closed smooth manifold for which there exists a measure $\mu$ compatible with inverse shadowing for $f$ and such that} $\mbox{Per}(f)\subset\mbox{supp}(\mu)$.

{\em Then the set} $\mbox{Int}^1({\mathcal I}_2)$ {\em coincides with the set of $\Omega$-stable difeomorphisms.}
\medskip

Indeed, by Lemma~4.3, if $f\in\mbox{Int}^1({\mathcal I}_2)$, then $f\in\mbox{Int}^1(\mbox{IS}(\mbox{Per}(f))$, and then $f$ is $\Omega$-stable (see the remark after Lemma~4.1).

On the other hand, if $f$ is $\Omega$-stable, then its set of periodic points is a dense subset of the hyperbolic set $\Omega(f)$; hence, $f\in\mbox{IS}(\mbox{Per}(f))$, and one may take as $\mu$ any measure whose support is $\Omega(f)$.
\medskip

Recall one more definition from Sec.~2.

For a homeomorphism $f$ of a metric space $K$, consider the set
$$\Psi(f)=\{p\in K:\, f\in\mbox{IIS}(\{p\})\}.$$ 

Thus, $f$ has the individual inverse shadowing property on a set $Y\subset K$ if and only if $Y\subset\Psi(f)$.

We say that a measure $\mu\in{\mathcal M}(K)$ is {\em compatible with individual inverse shadowing for $f$} if for any set $A\subset K$ with $\mu(A)>0$, 
$$A\cap\Psi(f)\neq\emptyset.$$

Denote by ${\mathcal I}_3$ the set of diffeomorphisms $f$ of a smooth closed manifold $M$ for which every $f$-invariant measure $\mu\in{\mathcal M}(M)$ is compatible with individual inverse shadowing.
\medskip

{\bf Theorem~4.2. }{\em The set} $\mbox{Int}^1({\mathcal I}_3)$ {\em coincides with the set of $\Omega$-stable diffeomorphisms.}
\smallskip

{\em Proof. } First we prove that any diffeomorphism $f\in\mbox{Int}^1({\mathcal I}_3)$ is $\Omega$-stable.

As above, it is enough to show that any periodic point of $f$ is hyperbolic. To get a contradiction, assume (as in the proof of Lemma~4.1) that we can find a fixed point $p$ of $f$, a diffeomorphism $h\in {\mathcal I}_3$, and a neighbourhood $U$ of $p$ with coordinates $(u,v)$ such that in $U$, $h$ has the form (4.2) in which the matrix $A$ is hyperbolic.

Take $a>0$ such that inclusion (4.1) is valid and fix $b,\ep_0>0$ such that
$$b+\ep_0<a.$$

Since the matrix $A$ is hyperbolic, there exists a neighbourhood $R_0$ of $p$ with the following properties: $R_0\subset R$ and if $h^k(q)=(s_k,t_k)$, then
$$|s_k|=|s_0|\leq b\mbox{ and }|t_k|\leq b,\quad q\in R_0,\eqno (4.9)$$
either for $k\geq 0$ or for $k\leq 0$.

Let us show that
$$R_0\cap\Psi(h)=\emptyset.\eqno (4.10)$$

Fix a point $q\in R_0$, assume that $q\in\Psi(h)$, fix an $\ep\in(0,\ep_0)$, and find the corresponding $d$. Find a $\de>0$ and a mapping $g\in C(M,M)$ such that (4.3) and (4.4) are satisfied. Then $\xi=\{g_k\}$, where $g_k=g,k\in\Z,$ is a $d$-method for $h$. Let $\{x_k=(u_k,v_k)\}$ be an arbitrary trajectory of the method $\xi$.

Clearly, there exist indices $k_1<0$ and $k_2>0$ such that for both $k=k_1$ and $k=k_2$:

-- either the point $x_{k}$ is outside $U$,

-- or $x_{k}\in U$ and $|u_{k}|\geq a$.

If relations (4.9) hold for $k\geq 0$, then
$$\dist(h^{k_1}(q),x_{k_1})>\ep;$$
otherwise,
$$\dist(h^{k_2}(q),x_{k_2})>\ep.$$
The contradiction obtained proves (4.10).

Now let us construct the required invariant measure. Consider the one-dimensional segment
$${\mathcal C}:=\{(u,v):\, u\in[-a,a],v=0\}.$$
Let $\mbox{mes}$ be the standard one-dimensional Lebesgue measure on ${\mathcal C}$ (so that  $\mbox{mes}({\mathcal C})=2a$). For an arbitrary measurable set $A\subset M$, set
$$\mu(A)=\frac{1}{2a}\mbox{mes}(A\cap{\mathcal C}).$$
Clearly, $\mu\in{\mathcal M}(M)$. Since every point of ${\mathcal C}$ is a fixed point of $h$, $h^{-1}(A\cap{\mathcal C})=A\cap{\mathcal C}$ for any set $A$; hence, the measure $\mu$ is $h$-invariant.

Since $R_0$ is an open set,
$$ \mu(R_0)=\frac{1}{2a}\mbox{mes}(R_0\cap{\mathcal C})>0,$$
and we get a contradiction between relation (4.10) and our assumption $h\in {\mathcal I}_3$. 

Thus, we have shown that any diffeomorphism  $f\in\mbox{Int}^1({\mathcal I}_3)$ is $\Omega$-stable.

Now let $f$ be an $\Omega$-stable diffeomorphism; denote by $\Omega(f)$ its nonwandering set. It is well known that $\mu(\Omega(f))=1$ for any $f$-invariant measure $\mu\in{\mathcal M}(M)$.

Hence, for any set $A$ with $\mu(A)>0$,
$$A\cap\Omega(f)\neq\emptyset;$$
it remains to note that any point $p\in A\cap\Omega(f)$ is a point of a hyperbolic nonwandering trajectory of $f$, hence, $p\in\Psi(f)$.

This completes the proof of the theorem. $\Box$
\bigskip

{\bf Acknowledgment. } This work was partially supported by the RFBR grant 18-01-00230-a.

\end{document}